\documentclass[12pt]{article}
\usepackage{amssymb}
\usepackage{amsfonts}
\usepackage{amsmath}
\usepackage{amsbsy}
\newcommand{\be}{\begin{equation}}
\newcommand{\ee}{\end{equation}}
\newcommand{\bea}{\begin{eqnarray}}
\newcommand{\eea}{\end{eqnarray}}
\newcommand{\ba}{\begin{array}}
\newcommand{\ea}{\end{array}}

\newcommand{\bc}{\begin{center}}
\newcommand{\ec}{\end{center}}
\newcommand{\ben}{\begin{enumerate}}
\newcommand{\een}{\end{enumerate}}
\newcommand{\bfi}{\begin{figure}}
\newcommand{\efi}{\end{figure}}

\newcommand{\bq}{\begin{quote}}
\newcommand{\eq}{\end{quote}}
\newcommand{\bqu}{\begin{quotation}}
\newcommand{\equ}{\end{quotation}}
\newenvironment{emphit}{\begin{itemize}}{\end{itemize}}
\newcommand{\bemp}{\begin{emphit}}
\newcommand{\eemp}{\end{emphit}}

\newcommand{\bt}{\begin{tabular}}
\newcommand{\et}{\end{tabular}}

\newtheorem{myth}{Theorem}[section]

\begin{document}
\date{}
\title{Not Every Osborn Loop Is Universal
}
\author{T. G. Jaiy\'e\d ol\'a\thanks{All correspondence to be addressed to this author.} \\
Department of Mathematics,\\
Obafemi Awolowo University,\\
Ile Ife 220005, Nigeria.\\
jaiyeolatemitope@yahoo.com\\tjayeola@oauife.edu.ng \and
J. O. Ad\'en\'iran \\
Department of Mathematics,\\
University of Agriculture, \\
Abeokuta 110101, Nigeria.\\
ekenedilichineke@yahoo.com\\
adeniranoj@unaab.edu.ng} \maketitle
\begin{abstract}
Michael Kinyon's 2005 open problem, based on the universality of
Osborn loops is solved. It is shown that not every Osborn loop is
universal.
\end{abstract}

\section{Introduction}
\paragraph{}
The universality of Osborn loops was raised as an open problem by Michael Kinyon \cite{phd33} in 2005. Up to this present moment, the problem is still open and it is our aim in this study to lay it to rest in this note. It will be shown here that not every Osborn loop is universal.

A loop is called an Osborn loop if it obeys the identity
\begin{equation*}
x(yz\cdot x)=x(yx^\lambda \cdot x)\cdot zx.
\end{equation*}
An Osborn loop is said to be universal if all its loop isotopes are Osborn loops.
The most popularly known varieties of
Osborn loops are CC-loops, Moufang loops, VD-loops and universal weak inverse property loops.
All these four varieties of Osborn loops are universal and this is
what makes the problem interesting.

For a comprehensive and detailed introduction and literature review of the subject of Osborn loops and universality, readers should check Jaiy\'e\d ol\'a and Ad\'en\'iran \cite{phd195}. The notations and symbols used here are exactly those adopted in \cite{phd195}.

\section{Main Result}
\begin{myth}\label{1:10}
Not every Osborn loop is universal.
\end{myth}
{\bf Proof}\\
Huthnance \cite{phd44} gave an example of an Osborn loop as follows.

Let $H=\mathbb{Z}\times \mathbb{Z}\times \mathbb{Z}$. Define a
binary operation $\star$ on $H$ by :
\begin{gather*}
[2i,k,m]\star [2j,p,q] = [2i+2j,k+p-ij(2j-1),q+m-ij(2j-1)]\\
[2i+1,k,m]\star [2j,p,q] =
[2i+2j+1,k+p-ij(2j-1)-j^2+j,q+m-ij(2j-1)-j^2]\\
[2i,k,m]\star [2j+1,p,q] = [2i+2j+1,m+p-ij(2j+1),q+k-ij(2j+1)]\\
[2i+1,k,m]\star [2j+1,p,q] =
[2i+2j+2,m+p-ij(2j+1)-j^2+j,q+k-ij(2j+1)-j^2]
\end{gather*}
for all $i,j,k,m,p,q\in \mathbb{Z}$. Assuming that $(H,\star )$ is a
universal Osborn loop, then it should obey the identity $v\cdot
vv=v^\lambda\backslash v\cdot v$ in Lemma~3.12 of \cite{phd195}. Let
$v=[2i+1,k,m]$. Then, by direct computation, we have
\begin{gather*}
v\cdot
vv=[6i+3,m+2k-10i^3-12i^2-2i,2m+k-10i^3-12i^2-i-1]~\textrm{and}\\
v^\lambda\backslash v\cdot
v=[6i+3,m+2k-14i^3-18i^2-7i-1,2m+k-14i^3-16i^2-6i-1].
 \end{gather*}
So, $v\cdot vv\ne v^\lambda\backslash v\cdot v$. Thus, $(H,\star )$
is not a universal Osborn loop.

\section{Concluding Remarks and Future Studies}
Theorem~\ref{1:10} completely lays to rest the open problem of
whether or not Osborn loops are generally universal. Kinyon
\cite{phd33} went furthermore to ask the following question, if in
case not every Osborn loop is universal.
\begin{itemize}
\item Does there exist a proper Osborn loop with trivial nucleus?
\end{itemize}
The answer to this question is due for future study.


\begin{thebibliography}{99}
\bibitem{phd44} Huthnance Jr., E.D., A theory of
generalised Moufang loops, Ph.D. thesis, Georgia Institute of
Technology, 1968.
\bibitem{phd195} Jaiy\'e\d ol\'a, T.G., and Ad\'en\'iran, J.O., New identities in universal Osborn loops, {\it Quasigroups And Related Systems}, \textbf{17} (2009), to appear.
\bibitem{phd33} Kinyon, M.K., A survey of Osborn loops,
\textit{Milehigh conference on loops, quasigroups and non-associative
systems}, University of Denver, Denver, Colorado, 2005.
\end{thebibliography}
\end{document}